\documentclass[english]{amsart}
\usepackage{ae}
\usepackage{aecompl}
\usepackage[T1]{fontenc}
\usepackage[latin1]{inputenc}
\usepackage{amssymb}

\makeatletter
 \theoremstyle{plain}
\newtheorem{thm}{Theorem}[section]
  \theoremstyle{plain}
  \newtheorem{lem}[thm]{Lemma}
  \theoremstyle{plain}
  \newtheorem{cor}[thm]{Corollary}

\usepackage[all,dvips,arc,curve]{xy}

\usepackage{mathbbol}

\usepackage{ifthen}

\newboolean{Details}

\setboolean{Details}{false}

\makeatother
\begin{document}

\title{A Projective $C^{*}$-Algebra Related to $K$-Theory }

\author{Terry A. Loring}

\email{loring@math.unm.edu}

\address{Department of Mathematics and Statistics, University of New Mexico,
Albuquerque, NM 87131, USA.}


\urladdr{http://www.math.unm.edu/\textasciitilde{}loring}

\begin{abstract}
The $C^{*}$-algebra $q\mathbb{C}$ is the smallest of the $C^{*}$-algebras
$qA$ introduced by Cuntz \cite{Cuntz-new-look-KK} in the context
of $KK$-theory. An important property of $q\mathbb{C}$ is the natural
isomorphism\[
K_{0}(A)\cong\lim_{\rightarrow}\left[q\mathbb{C},\mathbf{M}_{n}(A)\right].\]
Our main result concerns the exponential (boundary) map from $K_{0}$
of a quotient $B$ to $K_{1}$ of an ideal $I.$ We show if a $K_{0}$
element is realized in $\hom(q\mathbb{C},B)$ then its boundary is
realized as a unitary in $\tilde{I}.$ The picture we obtain of the
exponential map is based on a projective $C^{*}$-algebra $\mathcal{P}$
that is universal for a set of relations slightly weaker than the relations
that define $q\mathbb{C}$. A new, shorter proof of the semiprojectivity
of $q\mathbb{C}$ is described. Smoothing questions related the relations
for $q\mathbb{C}$ are addressed.
\end{abstract}
\maketitle
\begin{center}
\ifthenelse{\boolean{Details}}{{\Large Fully Detailed Version}}{}
\par\end{center}

\section{Introduction \label{sec:Introduction}}

The simplest nonzero projective $C^{*}$-algebra is $C_{0}(0,1].$
A quotient of this is $\mathbb{C},$ the simplest nonzero semiprojective
$C^{*}$-algebra. The first is universal for the relation $0\leq x\leq1$
and the second for $p^{*}=p^{2}=p.$ When lifting a projection from
a quotient, one must either settle for a lift that is only a positive
element or confront some $K$-theoretical obstruction to finding a
lift that is a projection. We consider noncommutative analogs of these
two $C^{*}$-algebras.

We use $\tilde{A}$ to denote the unitization of $A,$ where a unit
$\mathbb{1}$ is to be added even in $1_{A}$ exists. For elements
$h,$ $x$ and $k$ of $A,$ we use the notation\begin{equation}
T(h,x,k)=\left[\begin{array}{cc}
\mathbb{1}-h & x^{*}\\
x & k\end{array}\right]\in\mathbf{M}_{2}(\tilde{A}).\label{eq:T matrix def}\end{equation}

We will show that there is a $C^{*}$-algebra $\mathcal{P}$ with
generators $h,$ $k$ and $x$ that are universal for the relations\begin{align*}
 & hk=0,\\
 & 0\leq T(h,x,k)\leq1.\end{align*}
Moreover, $\mathcal{P}$ is projective. This does not appear to be
a familiar $C^{*}$-algebra, but it has a familiar quotient. The relations\begin{align*}
 & hk=0,\\
 & T(h,x,k)^{*}=T(h,x,k)^{2}=T(h,x,k)\end{align*}
have as their universal $C^{*}$-algebra the semiprojective $C^{*}$-algebra\[
q\mathbb{C}=\left\{ f\in C_{0}\left((0,1],\mathbf{M}_{2}\right)\left|f(1)\mbox{ is diagonal}\right.\right\} .\]

A complicated proof of the semiprojectivity of $q\mathbb{C},$ was
given in \cite{Loring-qC}. Subsequent proofs found with Eilers and
Pederson in \cite{ELP-anticommutation} and \cite{Loring-lifting-perturbing}
worked in the context of noncommutative CW-complexes. Those proofs
did not utilize the fact that $q\mathbb{C}$ is similar to the noncommutative
Grassmannian $G_{2}^{\mbox{nc}},$ c.f. \cite{Brown-free_product}.
The proof here uses this connection.

The importance of $q\mathbb{C}$ to $K$-theory is illustrated by
the isomorphism\[
K_{0}(A)\cong\left[q\mathbb{C},A\otimes\mathcal{K}\right]\cong\lim_{\rightarrow}\left[q\mathbb{C},\mathbf{M}_{n}(A)\right].\]
For example, see \cite{Cuntz-new-look-KK} and \cite{DadarlatLoringUnsuspendedE}.

Our main result concerns the exponential (boundary) map from $K_{0}$
of a quotient $B$ to $K_{1}$ of an ideal $I.$ If we look at $K_{0}$
as\[
K_{0}(D)\cong\lim_{\rightarrow}\left[q\mathbb{C},\mathbf{M}_{n}(D)\right]\]
then given\[
0\rightarrow I\rightarrow A\rightarrow B\rightarrow0\]
we show that a $K_{0}$ element realized in $\hom(q\mathbb{C},B)$
has boundary in $K_{1}(I)$ that can be realized as a unitary in the
$\tilde{I}.$ 

In the final section we look further into methods for perturbing approximate
representations of the relations for $q\mathbb{C}$ into true representations,
but this time restricting ourselves to using only $C^{\infty}$-functional
calculus.

\begin{lem}
\label{lem: low-level relations}The $C^{*}$-algebra\[
q\mathbb{C}=\left\{ f\in C_{0}\left((0,1],\mathbf{M}_{2}\right)\left|f(1)\mbox{ is diagonal}\right.\right\} \]
is universal in the category of all $C^{*}$-algebras for generators
$h,$ $k$ and $x$ with relations\begin{align}
h^{*}h+x^{*}x & =h,\nonumber \\
k^{*}k+xx^{*} & =k,\nonumber \\
kx & =xh,\nonumber \\
hk & =0.\label{eq:relations qC}\end{align}
The concrete generators may be taken to be\[
h_{0}=t\otimes e_{11},\quad k_{0}=t\otimes e_{22},\quad x_{0}=\sqrt{t-t^{2}}\otimes e_{21}.\]

\end{lem}
\begin{proof}
This is almost identical to Proposition 2.1 in \cite{Loring-qC}.
To see these are equivalent, notice first that the top two relations
imply $h$ and $k$ are positive. Since $x^{*}x$ is positive, the
relation $x^{*}x=h-h^{2}$ implies $h\leq1.$ It also implies $\| x\|\leq\frac{1}{2}.$
Similarly $k\leq1.$
\end{proof}
\begin{lem}
\label{lem: high-level relations}The $C^{*}$-algebra $q\mathbb{C}$
is universal in the category of all $C^{*}$-algebras for generators
$h,k,x$ and relations\begin{align}
 & hk=0,\nonumber \\
 & T(h,x,k)^{2}=T(h,x,k)^{*}=T(h,x,k).\label{eq:high-level-relations}\end{align}

\end{lem}
\begin{proof}
Since\[
T(h,x,k)=\left[\begin{array}{cc}
\mathbb{1}-h & x^{*}\\
x & k\end{array}\right]\]
and\[
T(h,x,k)^{2}=\left[\begin{array}{cc}
\mathbb{1}-2h+h^{2}+x^{*}x & x^{*}-hx^{*}+x^{*}k\\
x-xh+kx & k^{2}+xx^{*}\end{array}\right],\]
if we add $hk=0$ we have a set of relations equivalent to (\ref{eq:relations qC}).
\end{proof}

\section{Internal Matrix Structures in $C^{*}$-Algebras\label{sec:Internal-Matrix-Structures}}

\begin{lem}
\label{lem:internal-external matrices}Suppose $A$ is a $C^{*}$-algebra
and $X_{11},$ $X_{21},$ $X_{12},$ and $X_{22}$ are closed linear
subspaces of $A.$ Suppose $X_{ij}^{*}=X_{ji}$ and $X_{ij}X_{jk}\subseteq X_{ik}$
and $X_{11}X_{22}=0.$ 
\begin{enumerate}
\item The subset\[
\hat{X}=\left[\begin{array}{cc}
X_{11} & X_{12}\\
X_{21} & X_{22}\end{array}\right]\]
is a $C^{*}$-subalgebra of $\mathbf{M}_{2}(A).$
\item The sum\[
X_{11}+X_{21}+X_{12}+X_{22}\]
is a linear direct sum and is a $C^{*}$-subalgebra of $A,$ isomorphic
to $\hat{X}.$
\item There is a homotopy $\theta_{t}$ of injective $*$-homomorphisms
\[
\theta_{t}:X_{11}+X_{21}+X_{12}+X_{22}\rightarrow\mathbf{M}_{2}(A)\]
so that \[
\theta_{0}(x_{11}+x_{21}+x_{12}+x_{22})=\left[\begin{array}{cc}
x_{11}+x_{21}+x_{12}+x_{22} & 0\\
0 & 0\end{array}\right]\]
and\[
\theta_{1}(x_{11}+x_{21}+x_{12}+x_{22})=\left[\begin{array}{cc}
x_{11} & x_{12}\\
x_{21} & x_{22}\end{array}\right].\]

\end{enumerate}
\end{lem}
\begin{proof}
An element $x_{ij}$ of $X_{ij}$ factors as $x_{ij}=x_{ii}yx_{jj}$
with $y$ in $A$ and $x_{jj}=|x_{ij}^{*}|^{\frac{1}{4}}$ in $X_{jj}$
and $x_{ii}=|x_{ij}|^{\frac{1}{4}}$ in $X_{ii}.$ \ifthenelse{\boolean{Details}}{
If $x_{ij}$ is in $X_{ij}$ and $y_{kl}$ is in $X_{kl}$ with $j\neq k$
then for some $z$ and $w$ in $A$ and the factorization as above
we have \[
x_{ij}y_{kl}=x_{ii}zx_{jj}x_{kk}wx_{ll}=0\]
so $X_{ij}X_{kl}=0$ if $j\neq k.$ If $x$ is in $X_{ij}$ and in
$X_{kl}$ with $j\neq l$ then \[
xx^{*}\in X_{ij}X_{lk}=\{0\}\]
and so $xx^{*}=0$ and $x=0,$ proving $X_{ij}\cap X_{kl}=0$ when
$j\neq k.$ Similarly, this is true when $i\neq k.$}{From here,
it is easy to show that $X_{ij}X_{kl}=0$ if $j\neq k$ and that $X_{ij}\cap X_{kl}=0$
when $i\neq k$ or $j\neq l.$ $ $}

It is clear that $\hat{X}$ is a $C^{*}$-subalgebra of $\mathbf{M}_{2}(A).$
Let $ $$w_{t}$ be a partial isometry in $\mathbf{M}_{2}$ with $|w_{t}|=e_{11}$
for all $t$ and $w_{0}=e_{11}$ and $w_{1}=e_{21}.$ Define\[
\psi_{t}:\hat{X}\rightarrow A\otimes\mathbf{M}_{2}\]
 by\[
\psi_{t}\left(\sum x_{ij}\otimes e_{ij}\right)=\sum x_{ij}\otimes f_{ij}^{(t)}\]
 where\[
f_{11}^{(t)}=w_{t}^{*}w_{t},\quad f_{12}^{(t)}=w_{t}^{*}\]
 \[
f_{21}^{(t)}=w_{t},\quad f_{22}^{(t)}=w_{t}w_{t}^{*}.\]
 \ifthenelse{\boolean{Details}}{These $f_{ij}(t)$ satisfy\[
f_{ij}^{(t)}f_{jk}^{(t)}=f_{ik}^{(t)}\]
 but not the expected orthogonality relations. However, the orthogonality
of $X_{11}$ to $X_{22}$ and so forth gives us\[
\left(x_{ij}\otimes f_{ij}^{(t)}\right)\left(y_{kl}\otimes f_{kl}^{(t)}\right)=0\]
 and so\begin{eqnarray*}
\psi_{t}\left(\sum x_{ij}\otimes e_{ij}\right)\psi_{t}\left(\sum y_{ij}\otimes e_{ij}\right) & = & \left(\sum x_{ij}\otimes f_{ij}^{(t)}\right)\left(\sum y_{kl}\otimes f_{kl}^{(t)}\right)\\
 & = & \sum_{ijl}x_{ij}y_{jl}\otimes f_{il}^{(t)}\\
 & = & \sum_{il}\left(\sum_{j}x_{ij}y_{jl}\right)\otimes f_{il}^{(t)}\\
 & = & \psi_{t}\left(\sum_{il}\left(\sum_{j}x_{ij}y_{jl}\right)\otimes e_{il}\right)\\
 & = & \psi_{t}\left(\left(\sum x_{ij}\otimes e_{ij}\right)\left(\sum y_{kl}\otimes e_{kl}\right)\right)\end{eqnarray*}
}{The fact that $X_{ij}X_{kl}=0$ if $j\neq k$ implies that each
$\psi_{t}$ is a $*$-homomorphism. $ $}

The image of $\psi_{0}$ is \[
\left(X_{11}+X_{21}+X_{12}+X_{22}\right)\otimes e_{11}\]
and so we see that the direct sum of the $X_{ij}$ is a $C^{*}$-subalgebra
of $A.$

Now suppose\[
\psi_{t}\left(\sum x_{ij}\otimes e_{ij}\right)=0.\]
Then for all $r$ and all $s$ we have\[
0=\left(x_{rs}^{*}\otimes f_{1r}^{(t)}\right)\left(\psi_{t}\left(\sum_{ij}x_{ij}\otimes e_{ij}\right)\right)\left(x_{rs}^{*}\otimes f_{s1}^{(t)}\right)=x_{rs}^{*}x_{rs}x_{rs}^{*}\otimes e_{11}\]
which implies $x_{rs}=0.$ Therefore $\psi_{t}$ is injective.

If we let $\gamma$ denote the obvious isomorphism \[
\gamma:X_{11}+X_{21}+X_{12}+X_{22}\rightarrow\left(X_{11}+X_{21}+X_{12}+X_{22}\right)\otimes e_{11}\]
 and $\iota_{t}$ the inclusion of $\psi_{t}(\hat{X})$ into $\mathbf{M}_{2}(A)$
then\[
\theta_{t}=\iota_{t}\circ\psi_{t}\circ\psi_{0}^{-1}\circ\gamma\]
is \ifthenelse{\boolean{Details}}{a continuous path of injective
$*$-homomorphisms with \begin{eqnarray*}
\theta_{0}(x_{11}+x_{21}+x_{12}+x_{22}) & = & \iota_{t}\circ\gamma(x_{11}+x_{21}+x_{12}+x_{22})\\
 & = & (x_{11}+x_{21}+x_{12}+x_{22})\otimes e_{11}\end{eqnarray*}
and\begin{eqnarray*}
\theta_{1}(x_{11}+x_{21}+x_{12}+x_{22}) & = & \iota_{t}\circ\psi_{1}\circ\psi_{0}^{-1}\left((x_{11}+x_{21}+x_{12}+x_{22})\otimes e_{11}\right)\\
 & = & \iota_{t}\circ\psi_{1}\left(\sum_{ij}x_{ij}\otimes e_{ij}\right)\\
 & = & \sum_{ij}x_{ij}\otimes e_{ij}.\end{eqnarray*}
}{the desired path of injective $*$-homomorphisms.}
\end{proof}
\begin{lem}
\label{lem:linking alg plus units}Under the hypotheses of Lemma~\ref{lem:internal-external matrices},
the subset\[
\left[\begin{array}{cc}
\mathbb{C}\mathbb{1}+X_{11} & X_{12}\\
X_{21} & \mathbb{C}\mathbb{1}+X_{22}\end{array}\right]\]
is a $C^{*}$-subalgebra of $\mathbf{M}_{2}(\tilde{A}),$ and \[
\rho\left(\left[\begin{array}{cc}
\alpha\mathbb{1}+x_{11} & x_{12}\\
x_{21} & \alpha\mathbb{1}+x_{22}\end{array}\right]\right)=\alpha\oplus\beta\]
 determines a surjection onto $\mathbb{C}\oplus\mathbb{C}.$ 
\end{lem}
\begin{proof}
\ifthenelse{\boolean{Details}}{The sets\[
\left(\mathbb{C}\mathbb{1}+X_{11}\right)\otimes e_{11},\ X_{21}\otimes e_{21},\ X_{12}\otimes e_{12},\ \left(\mathbb{C}\mathbb{1}+X_{22}\right)\otimes e_{22}\ ,\]
subsets of $\mbox{M}_{2}(\tilde{A})$ satisfy the hypotheses of Lemma~\ref{lem:internal-external matrices},
and so their sum\[
\left[\begin{array}{cc}
\mathbb{C}\mathbb{1}+X_{11} & X_{12}\\
X_{21} & \mathbb{C}\mathbb{1}+X_{22}\end{array}\right]\]
is a $C^{*}$-subalgebra. }{This is follows easily from Lemma~\ref{lem:internal-external matrices}.}
\end{proof}
\begin{lem}
Suppose $I$ is an ideal in the $C^{*}$-algebra $A$ and $h$ and
$k$ in $A$ are positive elements. Then\[
I\cap\overline{kAh}=\overline{kIh}\]

\end{lem}
\begin{proof}
\ifthenelse{\boolean{Details}}{First consider the case $h=k.$ Since
$I$ is an ideal\[
\overline{hIh}\subseteq I\]
 and since $I$ is a subset of $A$ we also have\[
\overline{hIh}\subseteq\overline{hAh}.\]
 We thus have the easy direction,\[
\overline{hIh}\subseteq I\cap\overline{hAh}.\]

Suppose\[
x\in I\cap\overline{hAh}.\]
 e know there are continuous functions $g_{n}$ from the unit interval
to the unit interval, fixing $0,$ so that $g_{n}(h)h$ is an approximate
unit for $\overline{hAh}.$ Thus\[
hg_{n}(h)xg_{n}(h)h\rightarrow x\]
 and since\[
g_{n}(h)xg_{n}(h)\in I\]
 we have $x\in\overline{hIh}.$

For the general case, apply the special case to \[
\left[\begin{array}{cc}
h & 0\\
0 & k\end{array}\right],\ \left[\begin{array}{cc}
I & I\\
I & I\end{array}\right],\ \left[\begin{array}{cc}
A & A\\
A & A\end{array}\right].\]
 We learn that\[
\left[\begin{array}{cc}
I & I\\
I & I\end{array}\right]\cap\left[\begin{array}{cc}
\overline{hAh} & \overline{hAk}\\
\overline{kAh} & \overline{kAk}\end{array}\right]=\left[\begin{array}{cc}
\overline{hIh} & \overline{hIk}\\
\overline{kIh} & \overline{kIk}\end{array}\right].\]
}{The special case where $h=k$ is routine, and the general case
follows via a $2$-by-$2$ matrix trick.}
\end{proof}

\section{The Exponential Map in $K$-Theory\label{sec:The-Exponential-Map}}

We chose $b$ as the canonical generator of $K_{0}(q\mathbb{C})=\mathbb{Z},$
where $b$ is formed as the class of the projection \[
P_{0}=T(h_{0},x_{0},k_{0})\]
 minus the class of $[\mathbb{1}].$ (See (\ref{eq:T matrix def}).)

\begin{thm}
\label{thm:exponential using qC}Suppose 
\[
\xymatrix{
0  \ar[r] & 
	I \ar@{^{(}->}[r] &
		A \ar[r] ^{\pi} &
			B   \ar[r] &
				0 
}
\]
is
a short exact sequence of $C^{*}$-algebras. If $x$ is any element
of $K_{0}(B)$ such that $x=\varphi_{*}(b)$ for some $*$-homomorphism
$\varphi:q\mathbb{C}\rightarrow B,$ then $\partial(x)=[u]$ in $K_{1}(I)$
for some unitary $u\in\tilde{I}.$
\end{thm}
\begin{proof}
Let\[
y_{0}=\sqrt{t^{\frac{1}{2}}-t^{\frac{3}{2}}}\otimes e_{21}\]
so that $y_{0}$ is a contraction and\begin{equation}
k_{0}^{\frac{1}{8}}y_{0}h^{\frac{1}{8}}=x_{0}.\label{eq: x factors in qC}\end{equation}

Orthogonal positive contractions lift to orthogonal positive contractions,
so we can find $h$ and $k$ in $A$ with $\pi(h)=\varphi(h_{0}),$
$\pi(k)=\varphi(k_{0})$ $ $and\begin{align*}
 & hh=0,\\
 & 0\leq h\leq1,\\
 & 0\leq k\leq1.\end{align*}
Now take any $y$ in $A$ with $\pi(y)=\varphi(y_{0})$ and let $x=k^{\frac{1}{8}}yh^{\frac{1}{8}}$
and\[
T=T(h,x,k).\]
Then $\pi(x)=\varphi(x_{0}),$ \begin{equation}
\tilde{\pi}^{(2)}(T)=\tilde{\varphi}^{(2)}\left(P_{0}\right),\label{eq:lift of phi of P0}\end{equation}
\begin{equation}
T\in\left[\begin{array}{cc}
\mathbb{C}\mathbb{1}+\overline{hAh} & \overline{hAk}\\
\overline{kAh} & \mathbb{C}\mathbb{1}+\overline{kAk}\end{array}\right],\label{eq: in linking algebra}\end{equation}
\begin{equation}
\rho(T)=1\oplus0,\label{eq: added units are 1 0}\end{equation}
and $T^{*}=T.$

Let \[
f(\lambda)=\max(\min(\lambda,1),0)\]
and let $T^{\prime}=f(T).$ Then equations (\ref{eq:lift of phi of P0}),
(\ref{eq: in linking algebra}) and (\ref{eq: added units are 1 0})
hold with $T^{\prime}$ replacing $T.$ This means\[
T^{\prime}=T(h^{\prime},x^{\prime},k^{\prime})\]
for some $h^{\prime},$ $k^{\prime}$ and $x^{\prime}$ in $A$ that
are lifts of $h,$ $k$ and $x,$ and that\begin{align}
 & h^{\prime}k^{\prime}=0,\label{eq: relations for P}\\
 & 0\leq T\leq1.\nonumber \end{align}
This is an interesting lifting result that we will return to below.
For now, we turn to the exponential map.

Clearly $\partial([\mathbb{1}])=0$ so we need only compute $\partial\circ\varphi_{*}[P_{0}].$
We have the lifts $T$ and $T^{\prime}.$ We prefer to work with $T^{\prime}.$
A unitary that represents this $K_{1}$ element is $U^{\prime}=e^{2\pi iT^{\prime}}.$
Since\[
\tilde{\pi}^{(2)}\left(U^{\prime}\right)=\tilde{\varphi}^{(2)}\left(e^{2\pi iP_{0}}\right)=\left[\begin{array}{cc}
\mathbb{1} & 0\\
0 & \mathbb{1}\end{array}\right]\]
we know that \[
U^{\prime}\in\left[\begin{array}{cc}
\mathbb{1} & 0\\
0 & \mathbb{1}\end{array}\right]+\left[\begin{array}{cc}
I & I\\
I & I\end{array}\right].\]
By (\ref{eq: in linking algebra}) we know\[
U^{\prime}\in\left[\begin{array}{cc}
\mathbb{C}\mathbb{1}+\overline{hAh} & \overline{hAk}\\
\overline{kAh} & \mathbb{C}\mathbb{1}+\overline{kAk}\end{array}\right].\]
 Putting these facts together we discover\[
U^{\prime}\in\left[\begin{array}{cc}
\mathbb{1} & 0\\
0 & \mathbb{1}\end{array}\right]+\left[\begin{array}{cc}
\overline{hIh} & \overline{hIk}\\
\overline{kIh} & \overline{kIk}\end{array}\right]\subseteq\left[\begin{array}{cc}
\overline{hIh} & \overline{hIk}\\
\overline{kIh} & \overline{kIk}\end{array}\right]^{\sim}.\]
By Lemma~\ref{lem:internal-external matrices}, there is a path of
unitaries in $\left(\mathbf{M}_{2}(I)\right)^{\sim}$ from \[
U^{\prime}=\left[\begin{array}{cc}
u_{11} & u_{21}\\
u_{12} & u_{22}\end{array}\right]\]
to \[
\left[\begin{array}{cc}
-\mathbb{1}+u_{11}+u_{12}+u_{21}+u_{22} & 0\\
0 & \mathbb{1}\end{array}\right].\]
Thus $\partial\circ\varphi_{*}(b)=\partial\circ\varphi_{*}(P_{0})$
is represented in $\tilde{I}$ by the unitary\[
u=-\mathbb{1}+u_{11}+u_{12}+u_{21}+u_{22}.\]

\end{proof}
\begin{thm}
\emph{(\cite[Theorem 3.9]{Loring-qC})} $q\mathbb{C}$ is semiprojective.
\end{thm}
\begin{proof}
The proof of Theorem~\ref{thm:exponential using qC} is easily modified
to give a new proof of this result. One needs to assume that $I$
is the closure of the increasing union of ideals in $A.$ After the
lift $T$ is obtained in $B/I_{1},$ one can replace $I_{1}$ by $I_{n}$
with there now being a hole in the spectrum of $T$ around $\frac{1}{2}.$
Replacing the role of $f$ by \begin{equation}
f_{\frac{1}{2}}(\lambda)=\left\{ \begin{array}{ccc}
0 & \mbox{if} & \lambda<\frac{1}{2}\\
1 & \mbox{if} & \lambda\geq\frac{1}{2}\end{array}\right.,\label{eq:def of f1/2}\end{equation}
and following the same construction, one finds $T^{\prime}$ that
is a projection. The components of $T^{\prime}$ then provide a lift
in $B/I_{n}$ that is a representation of the generators of $q\mathbb{C}.$
\end{proof}
\begin{cor}
There is a universal $C^{*}$-algebra $\mathcal{P}$ for generators
$h,$ $k$ and $x$ for which\begin{align*}
 & hk=0,\\
 & 0\leq T(h,x,k)\leq1.\end{align*}
The surjection $\theta:\mathcal{P}\rightarrow q\mathbb{C}$ that sends
generators to generators is projective.
\end{cor}
\begin{proof}
Once we show $\mathcal{P}$ exists, the proof of the projectivity
of $\theta$ is contained in the proof of Theorem~\ref{thm:exponential using qC}.

By \cite{Loring-lifting-perturbing} we need only show that these
relations are invariant with respect of inclusions, are natural, are
closed under products, and are represented by a list of zero elements.
(This last requirement was erroneously missing in \cite{Loring-lifting-perturbing}.)
See also \cite{Hadwin-Kaonga-Mathes}.\ifthenelse{\boolean{Details}}{

If $\iota:A\rightarrow B$ is an inclusion, elements $h$ and $k$
of $A$ do or don't have product zero in $A$ or in $B.$ In addition,
the $*$-homomorphism\[
\tilde{\iota}^{(2)}:\mathbf{M}_{2}(\tilde{A})\rightarrow\mathbf{M}_{2}(\tilde{B})\]
is an injection. An element of a $C^{*}$-algebra is a positive contraction
even when considered as an element of a $C^{*}$-subalgebra. These
relations hold for $h,$ $k$ and $x$ in $A$ if and only if they
hold for the same $h,$ $k$ and $x$ in $B.$

If $\varphi:A\rightarrow B$ is a $*$-homomorphism, and elements
$h$ and $k$ of $A$ have product zero in $A,$ then they have product
zero in $B.$ If\[
0\leq T(h,x,k)\leq1\]
then\[
0\leq\tilde{\varphi}^{(2)}\left(T(h,x,k)\right)=T(\varphi(h),\varphi(x),\varphi(k))\leq1.\]
These relations hold for $\varphi(h),$ $\varphi(k)$ and $\varphi(x)$
in $B$ if they hold for $h,$ $k$ and $x$ in $A.$

As to the zero representation, $0$ is orthogonal to $0$ and\[
0\leq T(0,0,0)\leq1.\]

Suppose $x_{\lambda},$ $h_{\lambda}$ and $k_{\lambda}$ are in $A_{\lambda}$
and\[
h_{\lambda}k_{\lambda}=0\]
and\begin{equation}
0\leq T(h_{\lambda},x_{\lambda},k_{\lambda})\leq1\label{eq:order-relation-each-lambda}\end{equation}
for all $\lambda$ in $\Lambda.$ Then the norms of the $h_{\lambda},$
$k_{\lambda}$ and $x_{\lambda}$ are all bounded by $1,$ and so
in particular $\left\langle h_{\lambda}\right\rangle ,$$\left\langle k_{\lambda}\right\rangle $
and $\left\langle x_{\lambda}\right\rangle $ are all element of the
product $C^{*}$-algebra $\prod B_{\lambda}.$ In\[
\prod_{\lambda\in\Lambda}B_{\lambda}\]
we have\[
\left\langle h_{\lambda}\right\rangle \left\langle k_{\lambda}\right\rangle =\left\langle h_{\lambda}k_{\lambda}\right\rangle =\left\langle 0\right\rangle =0.\]
There is an inclusion of $C^{*}$-algebras\[
\Gamma:\mathbf{M}_{2}\left(\left(\prod_{\lambda}B_{\lambda}\right)^{\sim}\right)\hookrightarrow\prod_{\lambda}\mathbf{M}_{2}\left(\left(B_{\lambda}\right)^{\sim}\right)\]
given by\[
\Gamma\left(\left[\begin{array}{cc}
\alpha\mathbb{1}+\left\langle a_{\lambda}\right\rangle _{\lambda} & \beta\mathbb{1}+\left\langle b_{\lambda}\right\rangle _{\lambda}\\
\gamma\mathbb{1}+\left\langle c_{\lambda}\right\rangle _{\lambda} & \delta\mathbb{1}+\left\langle d_{\lambda}\right\rangle _{\lambda}\end{array}\right]\right)=\left\langle \left[\begin{array}{cc}
\alpha\mathbb{1}+a_{\lambda} & \beta\mathbb{1}+b_{\lambda}\\
\gamma\mathbb{1}+c_{\lambda} & \delta\mathbb{1}+d_{\lambda}\end{array}\right]\right\rangle _{\lambda}.\]
By (\ref{eq:order-relation-each-lambda}) we have\[
0\leq\left\langle T(h_{\lambda},x_{\lambda},k_{\lambda})\right\rangle _{\lambda}\leq1\]
and since $\Gamma$ is an inclusion, \[
0\leq T\left(\left\langle h_{\lambda}\right\rangle \left\langle x_{\lambda}\right\rangle _{\lambda},\left\langle k_{\lambda}\right\rangle _{\lambda}\right)\leq1.\]
Therefore $\left\langle h_{\lambda}\right\rangle _{\lambda},$$\left\langle k_{\lambda}\right\rangle _{\lambda}$
and $\left\langle x_{\lambda}\right\rangle _{\lambda}$ satisfy both
desired relations.}{ Details are left to the reader.}
\end{proof}
\begin{thm}
The $C^{*}$-algebra $\mathcal{P}$ is projective.
\end{thm}
\begin{proof}
Since $t^{2}\leq t$ in $C_{0}((0,1]),$ the matrix $T=T(h,x,k)$
satisfies $T^{2}\leq T.$ \ifthenelse{\boolean{Details}}{Therefore\[
\left[\begin{array}{cc}
\mathbb{1} & 0\\
0 & 0\end{array}\right]\left[\begin{array}{cc}
\mathbb{1}-h & x^{*}\\
x & k\end{array}\right]^{2}\left[\begin{array}{cc}
\mathbb{1} & 0\\
0 & 0\end{array}\right]\leq\left[\begin{array}{cc}
\mathbb{1} & 0\\
0 & 0\end{array}\right]\left[\begin{array}{cc}
\mathbb{1}-h & x^{*}\\
x & k\end{array}\right]\left[\begin{array}{cc}
\mathbb{1} & 0\\
0 & 0\end{array}\right].\]
}{}From this we deduce $x^{*}x\leq h-h^{2}.$ Similarly, \ifthenelse{\boolean{Details}}{we
have \[
\left[\begin{array}{cc}
0 & 0\\
0 & \mathbb{1}\end{array}\right]\left[\begin{array}{cc}
\mathbb{1}-h & x^{*}\\
x & k\end{array}\right]^{2}\left[\begin{array}{cc}
0 & 0\\
0 & \mathbb{1}\end{array}\right]\leq\left[\begin{array}{cc}
0 & 0\\
0 & \mathbb{1}\end{array}\right]\left[\begin{array}{cc}
\mathbb{1}-h & x^{*}\\
x & k\end{array}\right]\left[\begin{array}{cc}
0 & 0\\
0 & \mathbb{1}\end{array}\right],\]
from which we deduce}{}$xx^{*}\leq k-k^{2}.$ By \cite[Lemma 2.2.4]{Loring-lifting-perturbing}
we can factor $x$ as $x=k^{\frac{1}{8}}yh^{\frac{1}{8}}$ for some
$y$ in $\mathcal{P}.$ The rest of the proof is identical to argument
between equations (\ref{eq: x factors in qC}) and (\ref{eq: relations for P}).
\end{proof}

\section{Relations\label{sec:Relations}}

In this section we briefly examine a class of relations somewhat more
complicated than $*$-polynomials. See \cite{Hadwin-Kaonga-Mathes,Loring-lifting-perturbing,PhillipsInverseLimits}
for different approaches to relations in $C^{*}$-algebras,

Consider sets of relations of the form\[
f(p(x_{1},\ldots,x_{n}))=0,\]
either where $p$ is a self-adjoint $*$-polynomial in $n$ noncommuting
variables with $p(\mathbf{0})=0$ and \[
f\in C_{0}(\mathbb{R}\setminus\{0\}),\]
or where $p$ is not necessarily self-adjoint, $p(\mathbf{0})=0$
and $f$ is analytic on the plane. The point to restricting to these
relations is that\[
f(p(x_{1},\ldots,x_{n}))\]
makes sense, no matter the norm of the $C^{*}$-elements $x_{j},$
and so\[
\left\Vert f(p(x_{1},\ldots,x_{n}))\right\Vert \leq\delta\]
is a common-sense way to define an approximate representation.

Certainly a set $\mathcal{R}$ of relations on $x_{1},\ldots,x_{n}$
of this restricted form is invariant with respect to inclusion, is
natural, and each is satisfied when all the indeterminants are set
to $0.$ Therefore, $\mathcal{R}$ will define a universal $C^{*}$-algebra
if and only if it bounded, meaning for all $j$ we have\[
\sup\left\{ \left\Vert \tilde{x}_{j}\right\Vert \left|\strut\tilde{x}_{1},\ldots,\tilde{x}_{n}\mbox{ is a representation of }\mathcal{R}\right.\right\} <\infty.\]

We will also need to use relations of the form\begin{equation}
g\left(q\left(f_{1}(p_{1}(x_{1},\ldots,x_{n})),\ldots,f_{m}(p_{m}(x_{1},\ldots,x_{n}))\right)\right)=0\label{eq:typeOfRelation}\end{equation}
where the $f_{k},$ $p_{k}$ and $g,$ $q$ are pairs of continuous
functions and $*$-polynomials subscribing to the above rule. In particular
this will allow us the relation\[
\left\Vert q\left(f_{1}(p_{1}(x_{1},\ldots,x_{n})),\ldots,f_{m}(p_{m}(x_{1},\ldots,x_{n}))\right)\right\Vert \leq C.\]
For any $n$-tuple of elements in a $C^{*}$-algebra $A$ we define
$r(x_{1},\ldots,x_{n}),$ again in $A,$ by\[
r(x_{1},\ldots,x_{n})=f\left(q\left(f_{1}(p_{1}(x_{1},\ldots,x_{n})),\ldots,f_{m}(p_{m}(x_{1},\ldots,x_{n}))\right)\right).\]
If $x_{1},\ldots,x_{n}$ are is a sub-$C^{*}$-algebra, then so is
$r(x_{1},\ldots,x_{n}).$ Thus we are justified in the notation $r$
instead of the more pedantic $r_{A}.$ Also $r$ is natural. It is
still the case that the universal $C^{*}$-algebra exists if and only
if the set of relations is bounded.

\begin{lem}
\label{lem:secondaryRelations}Suppose\[
r_{k}(x_{1},\ldots,x_{n})=0\]
for $k=1,\ldots,K$ form a bounded set of relations of the form (\ref{eq:typeOfRelation}).
Suppose\[
s(x_{1},\ldots,x_{n})=0\]
is a relation of the form (\ref{eq:typeOfRelation}) that holds true
in\[
\mathcal{U}=C^{*}\left\langle x_{1},\ldots,x_{n}\left|r_{k}(x_{1},\ldots,x_{n})=0\ (\forall k)\right.\right\rangle .\]
Then for every $\epsilon>0$ there is a $\delta>0$ so that if $y_{1},\ldots,y_{n}$
in a $C^{*}$-algebra $A$ satisfy \[
\left\Vert r_{k}(y_{1},\ldots,y_{n})\right\Vert \leq\delta\quad(\forall k)\]
then\[
\left\Vert s(y_{1},\ldots,y_{n})\right\Vert \leq\epsilon.\]

\end{lem}
\begin{proof}
\ifthenelse{\boolean{Details}}{If the conclusion is false then there
are $C^{*}$-algebras $A_{c}$ and elements\[
x_{c,1},\ldots,x_{c,n}\in A_{c}\]
for which\[
\left\Vert r_{k}(x_{c,1},\ldots,x_{c,n})\right\Vert \leq\frac{1}{c}\]
for all $k$ and an $\epsilon_{0}>0$ so that\[
\left\Vert s(x_{c,1},\ldots,x_{c,n})\right\Vert \geq\epsilon_{0}.\]
Consider the elements\[
\left\langle x_{c,1}\right\rangle _{c},\ldots,\left\langle x_{c,n}\right\rangle _{c}\]
of\[
\prod_{c=1}^{\infty}A_{c}.\]
Let $\pi$ denote the quotient map\[
\pi:\prod A_{c}\rightarrow\prod A_{c}\left/\bigoplus A_{c}\right.\]
The naturality of $r,$ applied to the coordinate maps \[
\prod A_{c}\rightarrow A_{c_{0}},\]
gives us \[
r_{k}\left(\left\langle x_{c,1}\right\rangle _{c},\ldots,\left\langle x_{c,n}\right\rangle _{c}\right)=\left\langle r_{k}(x_{c,1},\ldots,x_{c,n})\right\rangle _{c}.\]
These will be elements of $\bigoplus A_{c}$ so we obtain a $*$-homomorphism\[
\varphi:\mathcal{U}\rightarrow\prod A_{c}\left/\bigoplus A_{c}\right.\]
with\[
\varphi(x_{j})=\pi\left(\left\langle x_{c,1}\right\rangle _{c}\right).\]
Then\begin{align*}
\pi\left(\left\langle s(x_{c,1},\ldots,x_{c,n})\right\rangle _{c}\right) & =s\left(\pi\left(\left\langle x_{c,1}\right\rangle _{c}\right),\ldots,\pi\left(\left\langle x_{c,n}\right\rangle _{c}\right)\right)\\
 & =s\left(\varphi(x_{1}),\ldots,\varphi(x_{n})\right)\\
 & =\varphi(s(x_{1},\ldots x_{n}))\\
 & =0\end{align*}
and so\[
\left\langle s(x_{c,1},\ldots,x_{c,n})\right\rangle _{c}\in\bigoplus A_{c}.\]
By the definition of a direct sum, this means\[
\lim_{c\rightarrow\infty}\left\Vert s(x_{c,1},\ldots,x_{c,n})\right\Vert =0,\]
a contradiction.}{This follows from standard arguments involving
the quotient of an infinite direct product by an infinite direct sum.}
\end{proof}

\section{Smoothing Relations}

We now modify the techniques from Section \ref{sec:The-Exponential-Map}
for a smooth version of semiprojectivity for $q\mathbb{C}.$ The result
is slightly weaker than \cite[Theorem 1.10]{Loring-qC}, but comes
with a more reasonable proof. The result involves maps from the generators
of $q\mathbb{C}$ to a dense $*$-subalgebra $A_{\infty}$ of a $C^{*}$-algebra
$A.$ The additional hypothesis is that $\mathbf{M_{2}}(A_{\infty}),$
and not just $A_{\infty},$ is closed under $C^{\infty}$ functional
calculus on self-adjoint elements. This additional assumption may
be no difficulty in examples. The smooth algebras of Blackadar and
Cuntz are closed under passing to matrix algebra (\cite[Proposition 6.7]{BlackadarCuntz}).

\begin{lem}
If $p^{*}=p$ is an element of a $C^{*}$-algebra $A$ and \[
\| p^{2}-p\|=\eta<\frac{1}{4}\]
then, with $f_{\frac{1}{2}}$ as in (\ref{eq:def of f1/2}), $f_{\frac{1}{2}}(p)$
is a projection in $A$ and\[
\left\Vert f_{\frac{1}{2}}(p)-p\right\Vert \leq\eta.\]

\end{lem}
\begin{proof}
This is well-known.
\end{proof}
\begin{thm}
For every $\epsilon>0,$ there is a $\delta>0$ so that if $A_{\infty}$
is a dense $*$-subalgebra of a $C^{*}$-algebra $A$ for which both
$A_{\infty}$ and $\mathbf{M_{2}}(A_{\infty})$ are closed under $C^{\infty}$
functional calculus on self-adjoint elements, then for any $h,$ $k$
and $x$ in $A_{\infty}$ for which\begin{align*}
 & \left\Vert h^{*}h+x^{*}x-h\right\Vert \leq\delta,\\
 & \left\Vert k^{*}k+xx^{*}-k\right\Vert \leq\delta,\\
 & \left\Vert kx-xh\right\Vert \leq\delta,\\
 & \left\Vert hk\right\Vert \leq\delta,\end{align*}
 there exist $\overline{h}$ $\overline{k}$ and $\overline{x}$ in
$A_{\infty}$ so that\begin{align*}
 & \overline{h}^{*}\overline{h}+\overline{x}^{*}\overline{x}-\overline{h}=0,\\
 & \overline{k}^{*}\overline{k}+\overline{x}\,\overline{x}^{*}-\overline{k}=0,\\
 & \overline{k}\,\overline{x}-\overline{x}\,\overline{h}=0,\\
 & \overline{h}\,\overline{k}=0,\end{align*}
and\[
\left\Vert \overline{h}-h\right\Vert \leq\epsilon,\quad\left\Vert \overline{k}-k\right\Vert \leq\epsilon,\quad\left\Vert \overline{x}-x\right\Vert \leq\epsilon.\]

\end{thm}
\begin{proof}
Let $\epsilon$ be given, with $0<\epsilon<\frac{1}{4}.$ Choose $\theta>0$
so that \[
\left\Vert h^{\prime}-h^{\prime\prime}\right\Vert \leq\theta,\quad\left\Vert k^{\prime}-k^{\prime\prime}\right\Vert \leq\theta,\quad\left\Vert x^{\prime}-x^{\prime\prime}\right\Vert \leq\theta,\]
\[
\left\Vert h^{\prime}\right\Vert \leq2,\quad\left\Vert k^{\prime}\right\Vert \leq2,\quad\left\Vert x^{\prime}\right\Vert \leq2,\]
implies\[
\left\Vert \left(h^{\prime*}h^{\prime}+x^{\prime*}x^{\prime}-h^{\prime}\right)-\left(h^{\prime\prime*}h^{\prime\prime}+x^{\prime\prime*}x^{\prime\prime}-h^{\prime\prime}\right)\right\Vert \leq\frac{\epsilon}{8},\]
\[
\left\Vert \left(k^{\prime*}k^{\prime}+x^{\prime}x^{\prime*}-k^{\prime}\right)-\left(k^{\prime\prime*}k^{\prime\prime}+x^{\prime\prime}x^{\prime\prime*}-k^{\prime\prime}\right)\right\Vert \leq\frac{\epsilon}{8},\]
\[
\left\Vert \left(k^{\prime}x^{\prime}-x^{\prime}h^{\prime}\right)-\left(k^{\prime\prime}x^{\prime\prime}-x^{\prime\prime}h^{\prime\prime}\right)\right\Vert \leq\frac{\epsilon}{8},\]
Choose $g_{+}$ some real-valued $C^{\infty}$ function on $\mathbb{R}$
for which\[
t\leq0\implies g_{+}(t)=0,\]
\[
t\geq0\implies t-\frac{\theta}{2}\leq g_{+}(t)\leq t,\]
and let $g_{-}(t)=g_{+}(-t).$ Choose $q_{+}$ some real-valued $C^{\infty}$
functions on $\mathbb{R}$ for which\[
t\leq0\implies q_{+}(t)=0,\]
\[
t\geq0\implies\sqrt{t-t^{2}}-\frac{\theta}{2}\leq\left(q_{+}(t)\right)^{2}\sqrt{t-t^{2}}\leq\sqrt{t-t^{2}},\]
and let $q_{-}(t)=q_{+}(-t).$ \ifthenelse{\boolean{Details}}{(Just
pick $q_{+}$ between $0$ and $1$ that is zero on the negative reals
and equal to $1$ on the positive reals above some small, positive
numbers.)}{}

Inside $q\mathbb{C},$ let $ $we have\[
g_{+}\left(\frac{1}{2}(h_{0}+h_{0}^{*}-k_{0}-k_{0}^{*})\right)=g_{+}(t)\otimes e_{11},\]
and\[
g_{-}\left(\frac{1}{2}(h_{0}+h_{0}^{*}-k_{0}-k_{0}^{*})\right)=g_{+}(t)\otimes e_{22}\]
and\begin{align*}
 & q_{-}\left(\frac{1}{2}(h_{0}+h_{0}^{*}-k_{0}-k_{0}^{*})\right)x_{0}q_{+}\left(\frac{1}{2}(h_{0}+h_{0}^{*}-k_{0}-k_{0}^{*})\right)\\
 & =\left(q_{+}(t)\right)^{2}\sqrt{t-t^{2}}\otimes e_{21}.\end{align*}
Therefore\[
\left\Vert g_{+}\left(\frac{1}{2}(h_{0}+h_{0}^{*}-k_{0}-k_{0}^{*})\right)-h_{0}\right\Vert \leq\frac{\theta}{2},\]
\[
\left\Vert g_{-}\left(\frac{1}{2}(h_{0}+h_{0}^{*}-k_{0}-k_{0}^{*})\right)-k_{0}\right\Vert \leq\frac{\theta}{2},\]
\[
\left\Vert q_{-}\left(\frac{1}{2}(h_{0}+h_{0}^{*}-k_{0}-k_{0}^{*})\right)x_{0}q_{+}\left(\frac{1}{2}(h_{0}+h_{0}^{*}-k_{0}-k_{0}^{*})\right)-x_{0}\right\Vert \leq\frac{\theta}{2}.\]
Of course, we also know \[
\| h_{0}\|\leq1,\quad\| k_{0}\|\leq1,\quad\| x_{0}\|\leq\frac{1}{2},\]
Lemma \ref{lem:secondaryRelations} tells us there is a $\delta>0$
so that if $h,$ $k$ and $x$ are in a $C^{*}$-algebra $A$ with\begin{align*}
 & \left\Vert h^{*}h+x^{*}x-h\right\Vert \leq\delta,\\
 & \left\Vert k^{*}k+xx^{*}-k\right\Vert \leq\delta,\\
 & \left\Vert kx-xh\right\Vert \leq\delta,\\
 & \left\Vert hk\right\Vert \leq\delta\end{align*}
then\[
\left\Vert g_{+}\left(\frac{1}{2}(h+h^{*}-k-k^{*})\right)-h\right\Vert \leq\theta,\]
\[
\left\Vert g_{-}\left(\frac{1}{2}(h+h^{*}-k-k^{*})\right)-k\right\Vert \leq\theta,\]
\[
\left\Vert q_{-}\left(\frac{1}{2}(h+h^{*}-k-k^{*})\right)xq_{+}\left(\frac{1}{2}(h+h^{*}-k-k^{*})\right)-x\right\Vert \leq\theta,\]
\[
\| h\|\leq2,\quad\| k\|\leq2,\quad\| x\|\leq2.\]
If necessary, replace $\delta$ with a smaller number to ensure $\delta<\frac{\epsilon}{2}.$

Let\[
\tilde{h}=f_{+}\left(\frac{1}{2}(h+h^{*}-k-k^{*})\right),\]
\[
\tilde{k}=f_{-}\left(\frac{1}{2}(h+h^{*}-k-k^{*})\right),\]
\[
h_{2}=g_{+}\left(\frac{1}{2}(h+h^{*}-k-k^{*})\right),\]
\[
k_{2}=g_{-}\left(\frac{1}{2}(h+h^{*}-k-k^{*})\right),\]
and\[
x_{2}=q_{-}\left(\frac{1}{2}(h+h^{*}-k-k^{*})\right)xq_{+}\left(\frac{1}{2}(h+h^{*}-k-k^{*})\right).\]
First notice that $\tilde{h}$ and $\tilde{k}$ are orthogonal positive
element of $A.$ Since\[
q_{+}\left(\frac{1}{2}(h+h^{*}-k-k^{*})\right)\]
is in the $C^{*}$-algebra generated by $\tilde{h},$ and\[
q_{-}\left(\frac{1}{2}(h+h^{*}-k-k^{*})\right)\]
is in the $C^{*}$-algebra generated by $\tilde{k},$ we have $x_{2}\in\overline{\tilde{k}A\tilde{h}}.$
Similarly, $h_{2}\in\overline{\tilde{k}A\tilde{h}}$ and $k_{2}\in\overline{\tilde{k}A\tilde{k}}.$
Next, observe that $h_{2},$ $k_{2}$ and $x_{2}$ are in $A_{\infty},$
with $h_{2}$ and $k_{2}$ self-adjoint and\[
\| h_{2}-h\|,\| k_{2}-k\|,\| x_{2}-x\|\leq\theta.\]
Therefore

\[
\left\Vert \left(h_{2}^{*}h_{2}+x_{2}^{*}x_{2}-h_{2}\right)-\left(h^{*}h+x^{*}x-h\right)\right\Vert \leq\frac{\epsilon}{8},\]
\[
\left\Vert \left(k_{2}k_{2}^{*}+x_{2}x_{2}^{*}-k_{2}\right)-\left(kk^{*}+xx^{*}-k\right)\right\Vert \leq\frac{\epsilon}{8},\]
\[
\left\Vert \left(k_{2}x_{2}-x_{2}h_{2}\right)-\left(kx-xh\right)\right\Vert \leq\frac{\epsilon}{8}\]
and so\[
\left\Vert h_{2}^{2}+x_{2}^{*}x_{2}-h_{2}\right\Vert \leq\delta+\frac{\epsilon}{8}\leq\frac{\epsilon}{4},\]
\[
\left\Vert k_{2}^{2}+x_{2}x_{2}^{*}-k_{2}\right\Vert \leq\delta+\frac{\epsilon}{8}\leq\frac{\epsilon}{4},\]
\[
\left\Vert k_{2}x_{2}-x_{2}h_{2}\right\Vert \leq\delta+\frac{\epsilon}{8}\leq\frac{\epsilon}{4}.\]

Let\[
T_{2}=T(h_{2},x_{2},k_{2})\in\left[\begin{array}{cc}
\mathbb{C}\mathbb{1}+\overline{\tilde{h}A\tilde{h}} & \overline{\tilde{h}A\tilde{k}}\\
\overline{\tilde{k}A\tilde{h}} & \overline{\tilde{k}A\tilde{k}}\end{array}\right].\]
With $\rho$ as in Lemma \ref{lem:linking alg plus units} $\rho\left(T_{2}\right)=1\oplus0.$
Since\[
\| T_{2}^{2}-T_{2}\|=\left\Vert \left[\begin{array}{cc}
-h_{2}+h_{2}^{2}+x_{2}^{*}x_{2} & x_{2}^{*}k_{2}-h_{2}x_{2}^{*}\\
k_{2}x_{2}-x_{2}h_{2} & -k_{2}+k_{2}^{2}+x_{x}x_{2}^{*}\end{array}\right]\right\Vert \]
we have\[
\| T_{2}^{2}-T_{2}\|\leq\frac{\epsilon}{2}.\]

Let $P=f_{\frac{1}{2}}(T_{2})$ and define $\overline{h},$ $\overline{k}$
and $\overline{x}$ via $T(\overline{h},\overline{x},\overline{k})=P.$
As in the proof of Theorem \ref{thm:exponential using qC} we see
that $x_{3},$ $k_{3}$ and $x_{3}$ satisfy the relations for $q\mathbb{C}.$
Since $f_{\frac{1}{2}}$ is smooth on intervals containing the spectrum
of $T_{2},$ these are elements of $A_{\infty}.$
\end{proof}

\end{document}